\documentclass[10pt, final]{amsproc}
\usepackage[cp850]{inputenc}
\usepackage{amssymb}
\usepackage[mathscr]{eucal}
\usepackage{amscd}
\usepackage{amsmath}
\usepackage{setspace}
\usepackage{tikz}
\theoremstyle{plain}

\newtheorem{prop}{Proposition}

\newtheorem{cor}{Corollary}

\theoremstyle{remark}

\author{Jes\'us Su\'arez de la Fuente}
\address{Departamento de Matem\'aticas, Avda. de Elvas s/n, Badajoz, Spain.}
\email{jesus@unex.es}

\thanks{The author was supported by project MTM2016-76958-C2-1-P and project IB16056 from La Junta de Extremadura.}

\title{Centralizers on a super-reflexive Schatten ideal}

\begin{document} 
\begin{abstract} We give a simple proof that there is no strictly singular bicentralizer on a super-reflexive Schatten ideal. This result applies, in particular, to the $p$-Schatten class for $1<p<\infty$.
\end{abstract}
\maketitle
\section{Main} Let $H$ be a separable Hilbert space and denote by $\mathcal B(H)$ the algebra of all bounded linear operators on $H$. We will follow the deep work of Kalton \cite{ka}, the interested reader may consider also \cite{CMR, G}, to introduce the following definitions. Let $E$ be a K\"othe sequence space which is symmetric, then we define the corresponding Schatten ideal $\mathcal C_E$ to be the algebra of all operators $A:H\to H$ whose singular values $s_n(A)$ satisfy that $(s_n(A))\in E$. We endow $\mathcal C_E$ with the norm $\|A\|_{E}=\| (s_n(A))\|_E$. In this notation the usual $p$-Schatten class is written as $C_{\ell_p}$. Let us recall that a map $\Omega: \mathcal C_E \to \mathcal B(H)$ is a \textit{bicentralizer} provided for some constant $k=k(\Omega)$ we have $$\|\Omega(VAW)-V\Omega(A)W\|_E\leq k \|V\|\|A\|_E\|W\|,$$
for $A\in \mathcal C_E$ and $V,W\in \mathcal B(H)$. We will also say that $\Omega$ is \textit{strictly singular} if its restriction to every infinite-dimensional closed subspace is never bounded. 

Regarding the $L_p$-spaces, the following remarkable result for centralizers was proved by Cabello S\'anchez \cite[Proposition 1]{C}.
\begin{prop}\label{ca}There is no strictly singular centralizer on $L_p$ for $0<p<\infty$.
\end{prop}
The proof of the proposition above contains several simplifications that are of independent interest. However, the core of the argument is based on a result of Kalton contained in the aforementioned paper \cite{ka}, namely that every centralizer on $L_p$ arises as a derivation from some interpolation process for $p>1$. The same work of Kalton contains a similar statement for the Schatten class \cite[Theorem 8.3]{ka}. The aim of this short note is to observe that using this last result one may prove a claim equivalent to Proposition \ref{ca} but for the Schatten class $\mathcal C_{\ell_p}$; and even overstep the limits of the $p$-class. 
\begin{prop}\label{ss}
There is no strictly singular bicentralizer on a super-reflexive Schatten ideal.
\end{prop}
\begin{proof} Let $\Omega$ denote a bicentralizer on $\mathcal C_{E}$. Then there exists hermitian bicentralizers $\Omega_1, \Omega_2$ and a constant $K>0$ so that if $A\in \mathcal C_{E}$ then
\begin{equation}\label{red}
\|\Omega(A)-\Omega_1(A)-i\Omega_2(A)\|_E\leq K\|A\|_E,
\end{equation}
see \cite[Lemma 8.1]{ka}. Let us find a copy of $\ell_2$ where $\Omega_1$ and $\Omega_2$ are bounded and the same will follow, by (\ref{red}), for $\Omega$. Thus, pick $\Omega_1$, our first hermitian bicentralizer, and apply \cite[Theorem 8.3]{ka} to find symmetric K\"othe sequence spaces $E_0, E_1$ such that $(\mathcal C_{E_0},\mathcal C_{E_1})_{\frac{1}{2}}=\mathcal C_{\widetilde E}$ (where $\widetilde E$ and $E$ have equivalent norms) and $\Omega_1$ is \textit{equivalent} to the bicentralizer induced by the previous interpolation formula, say $\widetilde \Omega_1$. We write this last claim in a formal language as there is some constant $K_1$ such that
\begin{equation}\label{red2}
\|\Omega_1(A)-\widetilde\Omega_1(A)\|_E\leq K_1\|A\|_E,
\end{equation}
for $A\in \mathcal C_E$. Now, let us give the copy of $\ell_2$ that makes the trick:
Pick $(e_n)_{n=1}^{\infty}$ an orthonormal basis in $H$. For this basis consider the span  of the corresponding first column matrix, or in other words the natural embedding
$$J: \ell_2 \longrightarrow  \mathcal C_E$$ $$a=(a_n)_{n=1}^{\infty} \longrightarrow A=\sum_{n=1}^{\infty} a_n (\;\cdot \;| e_1) e_n.$$
Writing
$$A=\left ( \sum_{n=1}^{\infty} a_n e_n \right) (\;\cdot \;| e_1)=\|(a_n)\|_{\ell_2}\left ( \sum_{n=1}^{\infty} \frac{a_n}{\|(a_n)\|_{\ell_2}} e_n \right) (\;\cdot \;| e_1),$$ one may easily deduce that $s_1(A)=\|(a_n)\|_{\ell_2}$ and $s_n(A)=0$ for $n\geq 2$. Therefore, the operator $J$ is well defined and the image spans a copy of $\ell_2$. Observe that $E$ does not play any special role here so that the same holds for $E_0$ and $E_1$. Indeed,
\begin{equation}\label{norm}
\|(s_1(A),0,0,...)\|_{E}=\|(s_1(A),0,0,...)\|_{E_0}=\|(s_1(A),0,0,...)\|_{E_1}=s_1(A).
\end{equation}  Thus, for a fixed $A=J(a)$, defining the constant function $F_A(z)=A$ where $z$ runs on the complex unit strip, the previous comments show that $F_A$ takes values in $\mathcal C_{E_0} \cap \mathcal C_{E_1}$. This gives that $F_A$ is a well defined map in the space of functions of Calder\'on $\mathcal F(\mathcal C_{E_0}, \mathcal C_{E_1})$ (see \cite{KM} for a clear exposition). And moreover, by (\ref{norm}), $$\|F_A\|_{\mathcal F(\mathcal C_{E_0}, \mathcal C_{E_1})}=\|A\|_E.$$ Hence, $F_A$ is a $C$-bounded selection of $A$ for the quotient map $\delta_{1/2}: \mathcal F(\mathcal C_{E_0}, \mathcal C_{E_1})\to C_{\widetilde E}$, where $C$ depends only on the equivalence constant between $E$ and $\widetilde E$.
Since $F_A$ is constant, its derivative must be zero, that in particular means $\widetilde\Omega_1(A)=0$ (see \cite[Page 1159]{KM} for a quick description of this phenomenon) and the same must hold for every $A\in J(\ell_2)$. Thus by (\ref{red2}) one has that $\Omega_1$ is bounded when restricted to $J(\ell_2)$, that is $$\|\Omega_1(A)\|_E\leq K_1\|A\|_E,$$ for every $A\in J(\ell_2)$. The same arguments work to find a constant $K_2$ satisfying $\|\Omega_2(A)\|_E\leq K_2\|A\|_E$ for every $A\in J(\ell_2)$. Using $K_1,K_2$ and the inequality (\ref{red}) one may easily check that $\Omega$ is bounded on $J(\ell_2)$ what finishes the proof.
\end{proof}
\begin{cor}
There is no strictly singular bicentralizer on $C_{\ell_p}$ for $1<p<\infty$.
\end{cor}

\end{document}